\input amstex.tex

\input amsppt.sty

\TagsAsMath

\magnification=1200

\hsize=5.0in\vsize=7.0in

\hoffset=0.2in\voffset=0cm

\nonstopmode

\document

\def\R{ \Bbb R}

\def\la{\langle}
\def\ra{\rangle}

\input amstex.tex
\input amsppt.sty
\TagsAsMath \NoRunningHeads \magnification=1200
\hsize=5.0in\vsize=7.0in \hoffset=0.2in\voffset=0cm \nonstopmode

\document

\topmatter

\title{A revision of "On
asymptotic stability in energy space of  ground states of NLS in
1D"}
\endtitle

\author
Scipio Cuccagna
\endauthor

\address
DISMI University of Modena and Reggio Emilia, via Amendola 2,
Padiglione Morselli, Reggio Emilia 42100 Italy\endaddress \email
cuccagna.scipio\@unimore.it \endemail

\abstract This is a revision of the author's paper "On asymptotic
stability in energy space of  ground states of NLS in 1D" \cite{C3}.
We correct an error in Lemma 5.4 \cite{C3} and we simplify the
smoothing argument.

\endabstract

\endtopmatter

 \head \S 1 Introduction \endhead

We consider even solutions  of a NLS
$$iu_t +  u_{xx} +  \beta ( |u|^2 )u=0=0\, ,\, (t,x)
\in \Bbb R \times \Bbb R   .\tag 1.1$$ We  assume $\beta (t)$
smooth, with {\item {(H1)}} $\beta (0)=\beta '(0)=0$, $\beta\in
C^\infty(\R,\R)$; {\item {(H2)}}  there exists a $p\in(1,\infty)$
such that for every $k=0,1$,
$$\left| \frac{d^k}{dv^k}\beta(v^2)\right|\lesssim
|v|^{p-k-1} \quad\text{if $|v|\ge 1$};$$ {\item {(H3)}} there exists
an open interval $\Cal{O}$ such that $
    u _{xx}-\omega u+\beta(u^2)u=0
$ admits a $C^1$-family of ground states $\phi _ {\omega }(x)$ for
$\omega\in\Cal{O}$;

{\item {(H4)}} $ \frac d {d\omega } \| \phi _ {\omega
}\|^2_{L^2(\Bbb R )}>0$ for $\omega\in\Cal{O}$.

\noindent By \cite{ShS} the $ \omega \to \phi _\omega \in H^1(\Bbb R
)$   is $C^2$ and by \cite{We1,GSS1-2} (H4) yields orbital stability
of the ground state $e^{i\omega t} \phi _ {\omega } (x)$. Here we
investigate asymptotic stability. We need some additional
hypotheses.

{\item {(H5)}} For any   $x\in \Bbb  {R} $, $u_0(x)=u_0(-x)$. That
is, the initial data $u_0$ of (1.1) are even.

{\item {(H6)}} Let $H_\omega$ be the linearized operator around
$e^{it\omega}\phi_\omega$, see (1.3).  $H_\omega$ has a positive
simple eigenvalue $\lambda(\omega)$ for $\omega\in\Cal{O}$. There
exists an $N\in\Bbb N$ such that
$N\lambda(\omega)<\omega<(N+1)\lambda(\omega)$. {\item {(H7)}} The
Fermi Golden Rule (FGR) holds (see Hypothesis 4.2
  in Section 4).
{\item {(H8)}} The point spectrum of $H_\omega$ consists of $0$ and
$\pm\lambda(\omega)$.   The points $\pm \omega$ are not resonances.

\proclaim{Theorem 1.1}    Let $\omega_0\in\Cal{O}$ and
$\phi_{\omega_0}(x)$ be a ground state. Let $u(t,x)$ be a solution
of (1.1). Assume (H1)--(H8). Then, there exist an $\epsilon_0>0$ and
a $C>0$ such that if $
\inf_{\gamma\in[0,2\pi]}\|u_0-e^{i\gamma}\phi_ {\omega }
\|_{H^1}<\epsilon <\epsilon _0,$ then there exist $\omega_+\in
\Cal{O}$, $\theta\in C^1(\R;\Bbb R)$ and $  h_+\in
 H^1 $    with $\| h_+\| _{H^1} \le C\epsilon  $ such that
$$\aligned &
\lim_{t\to\infty}\|u(t,\cdot)-e^{i\theta(t)}\phi_{\omega_{+}}
-e^{it\partial _x^2}h_+\|_{H^1}=0 .
 \endaligned $$
\endproclaim
Theorem  1.1  is the one dimensional version of  Theorem 1.1
\cite{CM}, which is valid for dimensions $D\ge 3$. In \cite{CM}
there is also a version of the theorem with  (H8) replaced by a more
general hypothesis, with more than one positive eigenvalue allowed
(but then a more restrictive (FGR) hypothesis (H7) is required). A
similar result could be proved here, but we prefer to skip the
proof. We recall  that  results of the sort discussed here were
pioneered by Soffer \& Weinstein \cite{SW1}, see also \cite{PW},
followed by Buslaev \& Perelman \cite{BP1-2}, about 15 years ago. In
this decade these early works were followed by a number of results
\cite{ BS,C1-2,GNT,M1-2,P,RSS,SW2,TY1-3,Wd1}.  It was heuristically
understood that the rate of the leaking of energy from the so called
"internal modes" into radiation, is small and decreasing when $N$
increases, producing technical difficulties in the closure of the
nonlinear estimates. For this reason prior to Gang Zhou \& Sigal
\cite{GS1}, the literature treated only the case when $N=1$ in (H6).
\cite{GS1} sheds light for $N>1$, with the eigenvalue $\lambda
(\omega )$ possibly very close to 0. Here we strengthen  the result
in \cite{GS1} for $D=1$, in analogy to the way   \cite{CM}
strengthens      \cite{GS1} for dimensions $D\ge 3$. For a detailed
introduction to the problem of asymptotic stability we refer to
\cite{CM}.    There are three hypotheses in \cite{GS1} which we
relax here. First of all, the (FGR) hypothesis in \cite{GS1} is more
restrictive than (H7). Specifically, \cite{GS1} require a sign
assumption on  a coefficient of  a certain equation
 obtained during a normal forms expansion. In \cite{CM} and later in
 this paper, it is shown that it is enough to assume that the
 coefficient be nonzero, a generic condition, and then it is proved
 that it has the right sign. Second, \cite{GS1} deals with solutions
 whose initial datum $u_0(x)$ satisfies more stringent conditions
 than being of finite energy. Finally, in the 1D case, \cite{GS1}
 requires that $\beta (t)$ be very small near 0, specifically
   $|\beta (t)|\lesssim |t|^{3N+2}$ for $|t|\le 1$, which we ease
   considerably here, since we only need $|\beta (t)|\lesssim
   |t|^2.$ Notice that the symmetry restriction (H5) is only required to
   avoid moving ground states, and that if we add to (1.1) some
   spacial inhomogeneity, thus eliminating translation invariance, then
   (H5) is unnecessary. So in particular our result, dropping (H5),
   will apply to equations like in \cite{GS1} of the form $ iu_t +  u_{xx}+V(x)u+\beta (    |u|^2) u=0
   $ with $V(x)$ a short range real valued potential.
As remarked in \cite{CM}, our result is
   relevant also to  equations of the form
   $ iu_t +  u_{xx}+V(x)    |u|^4 u=0 $ in the cases treated  Fibich and Wang \cite{FW}
   where ground states are proved to be orbitally stable.

The  proof of Theorem 1.1 is inspired  by Mizumachi \cite{M1} and
its use of Kato smoothing for the linearization which,      given $
\sigma _1=\left[ \matrix  0 &
1  \\
1 & 0
 \endmatrix \right]   , $ $
\sigma _2=\left[ \matrix  0 &
i  \\
-i & 0
 \endmatrix \right]  , $ $
\sigma _3=\left[ \matrix  1 &
0  \\
0 & -1
 \endmatrix \right]  , $
   is defined by

$$ H_\omega =\sigma _3 \left [ -d^2/dx^2 + \omega   -  \beta
(\phi ^2 _{\omega  }) - \beta ^\prime (\phi ^2 _{\omega  })\phi ^2
_{\omega  } \right ] +i  \beta ^\prime (\phi ^2 _{\omega  })\phi ^2
_{\omega  }.\tag 1.2$$ We exploit plane waves expansions for $
H_\omega $ and dispersive estimates  for the group $e^{-itH_\omega
}P_c( \omega ) $ proved in \cite{KS,GS1}. We also improve the
 Strichartz estimates  proved in \cite{KS} by means of a $TT^*$
 argument similar to the flat case.

We end with some notation. We set $\langle x \rangle =\sqrt{1+x^2}$.
We set $\| u \| _{H^{k,\tau}}:=\| \langle x \rangle  ^\tau u \|
_{H^{k }} .$   We set $ \langle f ,
  g\rangle =\int {^tf(x)} {\overline{g(x)}} dx,$  with $f(x)$ and $g(x)$ column vectors, $^tA$ the transpose and
  $\overline{g}$ the complex conjugate of $g$.
Given $x\in \R$ set $x^+=x\vee 0$ and $x^-=(-x)\vee 0$.
$R_H(z)=(H-z)^{-1}$.   $W^{k,p}(\Bbb R)$ is the space of tempered
distributions $f(x)$ such that $  (1-\partial _x^2) ^{k/2}f\in
L^p(\Bbb R)$.

\head \S 2 Linearization, modulation and set up \endhead

We will use  the  following classical result,  \cite{We1,GSS1-2}:

\proclaim{Theorem 2.1}  Suppose that $e^{i\omega t} \phi _ {\omega }
(x)$ satisfies (H4).   Then $\exists \, \epsilon >0$ and a
$A_0(\omega )>0$ such that for any  $\| u(0,x) - \phi _ {\omega } \|
_{H^1}<\epsilon $ we have for the corresponding solution $ \inf  \{
\|  u(t,x) -e^{i \gamma }\phi _ {\omega } (x-x_0) \| _{H^1(x\in \Bbb
R )} : \gamma \in \Bbb R \, \& \, x_0 \, \in \Bbb R  \}  <
A_0(\omega ) \epsilon . $ \endproclaim This statement is stronger
than the one in
   \cite{We1,GSS1-2} since  we state a more
precise estimate for the  $ \delta (\epsilon ) $ than in these
papers. We sketch the  proof   in \S 9.   Now we review some well
known facts about the linearization at a ground state. We can write
the ansatz $  u(t,x) = e^{i \Theta (t)} (\phi _{\omega (t)} (x)+
r(t,x)) \, , \, \Theta (t)= \int _0^t\omega (s) ds +\gamma (t). $
Inserting the ansatz into the equation we get
$$\aligned &
  i r_t  =
 -  r _{xx}+\omega (t) r-
\beta ( \phi _{\omega (t)} ^2 )r -\beta ^\prime ( \phi _{\omega (t)}
^2 )\phi _{\omega (t)} ^2 r \\&-
 \beta ^\prime ( \phi _{\omega (t)} ^2 )
\phi _{\omega (t)} ^2  \overline{  r }+ \dot \gamma (t) \phi
_{\omega (t)} -i\dot \omega (t)
\partial _\omega \phi   _{\omega (t)}
+ \dot \gamma (t) r
 +
  O(r^2).\endaligned
$$
 We set
$^tR= (r,\bar r) $,  $^t\Phi = ( \phi _{\omega } , \phi _{\omega } )
$
 and we rewrite the above equation as
$$i  R _t =H _{\omega}   R +\sigma _3 \dot \gamma   R
+\sigma _3 \dot \gamma \Phi - i \dot \omega \partial _\omega \Phi
+O(R^2).\tag 2.1
$$
 Set
  $H_0(\omega )=\sigma _3(-d^2/dx^2 +\omega )$
and $V(\omega )=H_\omega -  H_0(\omega ).$ The essential spectrum is
$$\sigma _e =\sigma _e (H_\omega )
=\sigma _e (H_0(\omega ) ) =(-\infty , -\omega ] \cup [ \omega ,
+\infty ) .$$  0 is an isolated eigenvalue. Given an operator $L$ we
set $N_g(L)= \cup _{j\ge 1} N(L^j)$ and $N(L)=\ker L$. \cite{We2}
 implies that,
 if $\{ \cdot \}$ means span,
$N_g(H^\ast _\omega )=\{ \Phi , \sigma _3\partial _\omega  \Phi   \}
$.  $\lambda  (\omega )$ has corresponding
   real eigenvector  $\xi  (\omega )$, which can be normalized  so that $\langle \xi   , \sigma _3 \xi   \rangle
 =1$.    $\sigma _1\xi (\omega ) $ generates $N(H  _\omega
+\lambda (\omega ))$ .  The function  $(\omega , x )\in \Cal O
\times \Bbb R  \to \xi  (\omega , x)$ is $C^2$; $|\xi  (\omega , x)|
< c e^{-a|x|}$ for fixed $c>0$ and $a>0$ if $\omega \in K \subset
\Cal O$, $K$ compact. $\xi   (\omega , x)$ is even in $x$ since by
assumption we are restricting ourselves in the category of such
functions. We have  the
 $H _{\omega}$   invariant
Jordan block decomposition
$$\align & L^2= N_g(H _{\omega} ) \oplus \big (
\oplus _{j,\pm  }N(  H _{\omega}\mp \lambda  (\omega )) \big )
\oplus L_c^2(H _{\omega})=N_g(H _{\omega})\oplus
 N_g^\perp (H _{\omega}^\ast  )
\endalign
$$
where we set $L_c^2(H _{\omega} )=\left \{   N_g(H _{\omega} ^\ast )
\oplus \oplus _{ \pm  }N(H _{\omega} ^\ast \mp \lambda (\omega
))\right \} ^{\perp}.$ We can impose
$$    R (t) =(z  \xi + \bar z  \sigma _1 \xi ) + f(t)   \in
\big [ \sum _{ \pm  } N(H _{\omega (t)}\mp \lambda  (\omega
(t)))\big ] \oplus L_c^2(H _{\omega (t)})  .\tag 2.2  $$ The
following claim admits an elementary proof which we skip:

\proclaim{Lemma 2.2} There is a Taylor expansion at $R=0$ of the
nonlinearity $O(R^2)$ in (2.1) with $R_{m,n}(\omega  ,x) $ and
$A_{m,n}(\omega ,x ) $ real vectors  and matrices rapidly decreasing
in $x$: $ O(R^2)=$
$$\aligned &  \sum _{ 2\le  m+n \le 2N+1} R_{m,n}(\omega ) z^m  \bar z^n+
\sum _{1\le  m + n \le  N} z^m  \bar z^n A_{m,n}(\omega ) f+
O(f^2+|z|^{2N+2}) .\endaligned  $$
\endproclaim
In terms of the frame in (2.2) and the expansion in Lemma 2.2, (2.1)
becomes
$$\aligned &if_t=\left ( H _{\omega (t)}+\sigma _3 \dot \gamma \right )f + \sigma _3
\dot \gamma \Phi (\omega )- i \dot \omega \partial _\omega \Phi (t)
+
  (z \lambda  (\omega ) -i\dot z ) \xi  (\omega ) \\&
- (\bar z \lambda  (\omega )+i\dot {\bar z }) \sigma _1\xi (\omega )
  +\sigma _3 \dot \gamma (z  \xi + \bar z  \sigma _1 \xi )
-i \dot \omega (z  \partial _\omega \xi + \bar z  \sigma _1
\partial _\omega \xi ) \\& + \sum _{ 2\le  m+n \le 2N+1} z^m \bar
z^nR_{m,n}(\omega ) + \sum _{1\le  m + n \le N} z^m \bar z^n
A_{m,n}(\omega ) f+\\& + O(f^2)+ O_{loc}(|z  ^{2N +2}|)
\endaligned \tag 2.3
$$
where by $O_{loc}$ we mean that the there is a factor $\chi (x)$
rapidly decaying   to 0   as $|x|\to \infty $. By taking inner
product of the equation with generators of $N_g(H_\omega ^\ast )$
and $N(H_\omega ^\ast -\lambda  )$ we obtain modulation and discrete
modes equations:

$$\aligned & i\dot \omega \frac{d\| \phi _\omega \| _2^2}{d\omega  }
=\langle \sigma _3 \dot \gamma (z \xi + \bar z  \sigma _1 \xi ) -i
\dot \omega (z \partial _\omega \xi + \bar z  \sigma _1
\partial _\omega \xi )  +   \sum _{    m+n =2}^{ 2N+1}
 z^m \bar z^nR_{m,n}(\omega )\\& +  \big ( \sigma _3 \dot \gamma +i\dot \omega
\partial _\omega P_c+ \sum _{  m + n =1}^{ N} z^m \bar z^n
A_{m,n}(\omega )  \big ) f   + O(f^2)+ O_{loc}(|z ^{2N +2}|), \Phi
\rangle \\& \dot \gamma \frac{d\| \phi _\omega \| _2^2}{d\omega  }
=\langle \text{ same as above }, \sigma _3
\partial _\omega \Phi \rangle \\& i\dot z -\lambda (\omega )z =\langle
\text{ same as above }, \sigma _3 \xi  \rangle .
\endaligned \tag 2.4$$

\head \S 3 Spacetime estimates for $H_\omega $ \endhead

We collect some linear estimates needed for the proof of Theorem 1.1
in  \S 4. First of all we prove that the group $e^{-itH_{\omega
}}P_c(\omega )$ satisfies the same Strichartz estimates of the flat
case. The proof is almost the same of the flat case. In particular
we are able to implement a $TT^\ast $ argument. For a different
proof without the $L_t^4L_x^\infty$ estimate, see Corollary 7.3
\cite{KS}.

\proclaim{Lemma 3.1 (Strichartz estimate)} There exists a positive
number $C=C(\omega )$ upper semicontinuous in $\omega $ such that
for any $k\in [0,2]$:
  {\item {(a)}}
 for any $f\in
L^2_c( {\omega })$,
$$\|e^{-itH_{\omega }} f\|_{L_t^4W_x^{k,\infty }\cap L_t^\infty H_x^k }\le C\|f\|_{H^k}.$$
 {\item {(b)}}
  for any $g(t,x)\in
S(\R^2)$,
$$
\|\int_{0}^te^{-i(t-s)H_{\omega }} P_c( {\omega
})g(s,\cdot)ds\|_{L_t^4W_x^{k,\infty }\cap L_t^\infty H_x^k} \le
C\|g\|_{L_t^{4/3}W_x^{k,1 }+L_t^1H_x^k}.
$$\endproclaim
{\it Proof.} First of all, the case $0<k\le 2$ follows by the case
$k=0$ by a simple argument in Corollary 7.3 \cite{KS}. Now we focus
on the $k=0$ case. For any $2\le p \le \infty$ by
 \cite{BP1,KS,GS2}  $\exists$ $C=C(\omega )$ upper semicontinuous in
$\omega $ such that
$$\|e^{-itH_{\omega }}  P_c( {\omega
})f\|_{ L_x^p} \le C  t^{-\frac{1}{2}+ \frac{1}{p}}\|f\|_{L
^{\frac{p}{p-1}} }. \tag 1$$  $(b)$ is a consequence of (1) and of
Hardy Littlewood theorem. The $L_t^\infty L_x^2$ estimate in  $(a)$
is an immediate consequence of (1) for $p=2$. The quadratic form
$\langle f , \sigma _3g\rangle $ defined in $ L^1_c( {\omega
})\times L^\infty _c( {\omega } )$ establishes an isomorphism
$(L^1_c( {\omega })) ^{\ast} \simeq  L^\infty _c( {\omega } )$.
Based on  $\langle e^{-itH_{\omega }}f , \sigma _3 g\rangle =
\langle f , \sigma _3  e^{ itH_{\omega }}g \rangle $  the following
operators  are formally adjoints
$$ \aligned & g(t,x)\in L_t^{4/3}L_c^1( {\omega })\to T  g=\int _{\R} e^{itH_{\omega
}} g(t) (x) dt \in L^2_c( {\omega }) \\& \text{and } \, f\in L^2_c(
{\omega })\to T^\ast f=e^{-itH_{\omega }} f \in L_t^4L_c^\infty (
{\omega }).\endaligned $$ Then we can perform a slight modification
of the standard $TT^\ast$ argument. Preliminarily, we split $P_c(
\omega )= P_+( \omega )+P_-( \omega )$ the projections in the
positive and negative part of $\sigma _c(H_\omega )$, see Appendix B
 and \cite{BP2,BS,C2}. We bound separately $  P_\pm ( \omega
)\circ Tf .$     The operator $T^\ast \circ P_\pm ( \omega )\circ T$
  is   bounded   thanks to (1) and Hardy
Littlewood theorem. We write, for $L_c^p=L_c^p( \omega )$,
$$ \aligned & |\langle  P_\pm ( \omega )\circ T f,\sigma _3
 P_\pm ( \omega )\circ T f\rangle _{tx}|  = |\langle T^\ast \circ P_\pm ( \omega )\circ T
f,\sigma _3 f\rangle  _{tx}| \le \\& \le \| T^\ast \circ P_\pm (
\omega )\circ T : L_t^{4/3}L_c^1 \to L_t^4L_c^\infty \| \, \| f \|
_{L_t^{4/3}L_c^1} ^2 .\endaligned $$ Assuming
$$\langle  P_\pm ( \omega )h,\sigma _3
 P_\pm ( \omega )h\rangle _{ x} \approx \pm \|  P_\pm ( \omega )h \| _{L^2_c( \omega )} \tag 2$$
we conclude $  \|  P_\pm ( \omega )\circ T
 f \| _{L^2_c( \omega )} \lesssim \| f
\| _{L_t^{4/3}L_c^1( \omega )} .$
 Adding up  we get $\|    T
 f \| _{L^2_c( \omega )} \lesssim \| f
\| _{L_t^{4/3}L_c^1( \omega )}  $. For $\psi \in C_0([0,\infty
)\times \Bbb R)$ we get the following which yields (a):
$$\aligned & \langle   T ^\ast f,\sigma _3 \psi \rangle _{tx} = \langle     f,\sigma _3 T\psi \rangle _{tx}
\le C \| f\| _{L^2_c( \omega )} \| \psi \| _{ L_t^{4/3}L_c^1}.
\endaligned $$
 To
obtain (2) we observe that there exists a wave operator $ W
:L^2(\Bbb R)\to L^2_c( \omega )$ which is an isomorphism with
inverse $Z$ such that  for $h=W\widetilde{h}$ and
$^t\widetilde{h}=(\widetilde{h}_1,\widetilde{h}_2)$ we have
 $$\aligned & \langle  P_+ ( \omega )h,\sigma _3
 P_+ ( \omega )h\rangle =  \| \widetilde{h}_1 \| ^2 _2 \approx \|  P_+ ( \omega )h \| _{L^2_c( \omega )} \text{ and}\\&\langle  P_- ( \omega )h,\sigma _3
 P_- ( \omega )h\rangle = - \| \widetilde{h}_2 \| ^2 _2 \approx -\|  P_- ( \omega )h \| _{L^2_c( \omega )}.
 \endaligned$$
$W$   and $Z$ above can be defined in a standard way, $Z$ thanks to
(1) and Proposition 8.1 \cite{KS}, as strong limits $ W(\omega
)=\lim _{t\to +\infty} e^{-itH_\omega}e^{it\sigma _3(-\Delta +\omega
)}$, $ Z(\omega )=\lim _{t\to +\infty} e^{it\sigma _3( \Delta
-\omega )} e^{ itH_\omega }$ and by standard theory they are
inverses of each other.

\proclaim{Lemma 3.2} Fix $\tau >3/2$.

{\item {(1)}} There exists $C=C(\tau ,\omega )$, upper
semicontinuous in $\omega $ such that  for any $\varepsilon \neq 0$

$$\| R_{H_\omega }(\lambda +i\varepsilon )P_c(H_\omega )
u\| _{L^2_\lambda L^{2,-\tau }_x}\le  C \| u\| _{L^2}. $$ {\item
{(2)}} For any $u\in L^{2, \tau }_x $ the following limits:
$$ \lim _{\epsilon \searrow 0}R_{H_\omega }(\lambda \pm i\varepsilon )
u= R_{H_\omega }^\pm  (\lambda ) u  \text{ in $C^0(\sigma
_e(H_\omega ),L^{2, -\tau }_x)$}.$${\item {(3)}} We have
$$
  \|   R_{H_\omega }^\pm  (\lambda
)P_c(H_\omega )   \| _{B( L^{2,\tau }_x, L^{2,-\tau }_x)} < C
\langle \lambda \rangle ^{-\frac{1}{2}} .$$ {\item {(4)}} Given any
$u\in L^{2, \tau }_x $ we have
$$P_c(H_\omega )u=\frac{1}{2\pi i}\int _{\sigma _e(H_\omega )}
(R_{H_\omega }^{+}(\lambda  )-R_{H_\omega }^{-}(\lambda  ))  u\,
d\lambda .$$

\endproclaim
These are consequences of the fact that $\sigma _e(H_\omega )$ does
not contain eigenvalues and that $\pm \omega $ are not resonances,
and of the theory in \cite{KS}.

 \proclaim{Lemma 3.3} For any $k$ and  $\tau >3/2$  $\exists$
  $C=C(\tau ,k,\omega )$ upper
semicontinuous in $\omega $ such that:
   {\item {(a)}}
  for any $f\in S(\R)$,
$$\align &
\| e^{-itH_{\omega }}P_c(H_\omega )f\| _{L_{  t}^2 H_x^{k, -\tau}}
\le
 C\|f\|_{H ^{k}} .
\endalign $$
  {\item {(b)}}
  for any $g(t,x)\in
 {S}(\R^2)$
$$ \left\|\int_\Bbb R e^{itH_{\omega }}
P_c(H_\omega )g(t,\cdot)dt\right\|_{H^k_x} \le C\| g\|_{L_{  t}^2
H_x^{k, \tau}}.
$$\endproclaim
  {\it Proof.} It is enough to prove
  Lemma 3.3, as well as Lemmas 3.4 below, for $k=0$.
   (a) implies (b) by duality:  $$\aligned & |\langle  f, \sigma _3\int_\R e^{itH_{\omega }}P_c( {\omega })g(t )dt
\rangle _{x}|= |\langle \la x\ra^{-\tau}e^{-itH_{\omega }}P_c (
{\omega })f, \sigma _3\la x\ra^{ \tau}g \rangle _{tx}|\\& \le \|
e^{-itH_{\omega }}P_c( {\omega })f\|_{L_{  t}^2 L_x^{2, -\tau}}\|
g\|_{L_{  t}^2 L_x^{2, \tau}} \le \| f\|_{L_x^2}\| g\|_{L_{ t}^2
L_x^{2, \tau}}.
\endaligned
$$ We now prove (a) for $k=0$.
Let $g(t,x) \in S(\Bbb R^2) $ with $g(t)=P_c(H_\omega )g(t)$. Then
$$\aligned & \la  e^{-itH_\omega } f,\sigma _3
g\ra _{t,x}=  \frac{1}{\sqrt{2\pi}i}\int_\R e^{-i\lambda t}
  \left\la
 (R_{H_\omega }^{+}(\lambda  )-R_{H_\omega
}^{-}(\lambda  ))   f,\sigma _3 \overline{\widehat{g}}(\lambda
)\right\ra_x  d\lambda  \\& = \frac{1}{\sqrt{2\pi}i}\int_{ \sigma
_e(H_\omega )} e^{-i\lambda t}
  \left\la
 (R_{H_\omega }^{+}(\lambda  )-R_{H_\omega
}^{-}(\lambda  ))   f,\sigma _3 \overline{\widehat{g}}(\lambda
)\right\ra_x  d\lambda  .
\endaligned $$
Then from Fubini and  Plancherel and by   (1) Lemma 3.3 we have
$$\aligned & \big |\la  e^{-itH_\omega } f,\sigma _3 g\ra _{t,x}\big
|
 \le (2\pi)^{-1/2}
\| (R_{H_\omega }^{+}(\lambda  )-R_{H_\omega }^{-}(\lambda ))
f\|_{L^{2, -\tau}_x(\Bbb R)L^2_\lambda (\sigma _e(H_\omega ))}
\times \\& \times \| \widehat{g}(\lambda,\cdot) \|_{L^{2, \tau}_x
L^2_\lambda  }
  \lesssim
\|  f \|_{L_x^  2   } \|  {g} \|_{L^{2, \tau}_x L^2_t  }
 .
  \endaligned $$

\proclaim{Lemma 3.4} For any $k$ and  $\tau >3/2$  $\exists$
  $C=C(\tau ,k,\omega )$ as
above such that $\forall$ $g(t,x)\in {S}(\R^2)$
$$\align &  \left\|  \int_0^t e^{-i(t-s)H_{\omega
}}P_c(H_\omega )g(s,\cdot)ds\right\|_{L_{  t}^2 H_x^{k, -\tau}} \le
C\|  g\|_{L_{  t}^2 H_x^{k, \tau}}.\endalign
$$
\endproclaim
{\it Proof.}  By Plancherel and H\"older inequalities and by (3)
Lemma
 3.2   we have
$$ \aligned &
\| \int _{0}^t e^{-i(t-s)H_{\omega }}P_c(H_\omega )g(s,\cdot)ds\|_{
L_{t }^2L_{ x}^{2,-\tau }} \le \\& \le  \| R_{H_\omega }^+(\lambda
)P_c (H_\omega )
  \widehat{ \chi }_{[0,+\infty )}\ast _\lambda
   \widehat{ g}(\lambda,x)\|_{L_{t }^2L_{ x}^{2,-\tau } }   \le \\& \le
   \left\| \,
\|   R_{H_\omega }^+ (\lambda )P_c  (H_\omega ) \| _{B(
L^{2,\tau}_x, L^{2,- \tau}_x)} \|
   \widehat{ \chi }_{[0,+\infty )}
   \ast _{\lambda } \widehat{g} (\lambda,x) \|_{L_{ x}^{2, \tau }}\, \right\|_{L^2_\lambda}
\\ \le &
  \|  R_{H_\omega }^+ (\lambda
)P_c(H_\omega )   \| _{L^\infty _\lambda (\Bbb R ,B( L^{2,\tau}_x,
L^{2,-\tau}_x))}\| g\|_{L_{t }^2L_{ x}^{2, \tau } }  \le C \|
g\|_{L_{t }^2L_{ x}^{2, \tau } } .
\endaligned $$

\proclaim{Lemma 3.5} $k$ and  $\tau >3/2$  $\exists$
  $C=C(\tau ,k,\omega )$ as
above such that $\forall$ $g(t,x)\in {S}(\R^2)$
$$\align  &
\left\|\int_0^t e^{-i(t-s)H_{\omega }}P_c(H_\omega )g(s,\cdot)ds
  \right\|_{
L_t^\infty L_x^2\cap  L^{4 }_{t}(\Bbb R ,W ^{k,\infty }_x) } \le
C\|g\|_{L_t^2H_x^{k,\tau}}.
\endalign$$
\endproclaim
 {\it Proof.} For $g(t,x)\in S(\Bbb R^2)$ set
$$T g(t)=\int _0^{+\infty} e^{-i(t-s)H_{\omega }}
P_c(H_\omega )g(s) ds.$$ Lemma 3.3 (b) implies $f:=\int _0^{+\infty}
e^{isH_\omega }P_c( \omega )g(s)ds\in L^2(\Bbb R)$. Then Lemma 3.5
 is a direct consequence of \cite{CK}.

\head \S 4 Proof of Theorem 1.1\endhead

   We restate Theorem 1.1 in a more precise form:

\proclaim{Theorem 4.1} Under the assumptions of Theorem 1.1   we can
express
$$u(t,x)=e^{i \Theta (t)} \left (\phi _{\omega (t)} (x)+  \sum _{j=1}^{2N} p_j(z,\bar z)  A _j (x,\omega (t))+
 h(t,x)\right )$$ with $p_j(z,\bar z)=O(z )$ near 0,
with  $\lim _{t\to +\infty}\omega (t)$ convergent, with $ |A _j
(x,\omega (t)) |\le C   e^{-a|x|}$ for fixed $C>0$ and $a>0$, $\lim
_{t\to +\infty} z(t  )  =0,$ and   for fixed $C>0$
$$\| z(t)\|_{ L_t ^{2N+2}}^{N+1} +\| h(t,x)\| _
{L^\infty  _tH ^{1 }_x \cap L^5 _tW ^{1,10}_x \cap L^4_tL^\infty _x}
<C \epsilon .\tag 1$$
   Furthermore, there exists
$h_\infty\in H^1(\Bbb R ,\Bbb C )$ such that
$$\lim_{t\to\infty}\|  e^{i\int _0^t\omega (s) ds +i\gamma (t)}h(t)-
e^{it  \frac{d^2}{dx^2}   }h_\infty\|_{H^1}=0. \tag 2$$
\endproclaim

The proof of Theorem 4.1  consists in a normal forms expansion and
in the closure of some nonlinear estimates. The normal forms
expansion is exactly the same of \cite{CM}, in turn an adaptation of
\cite{GS1}.

\head \S 4.1 Normal form expansion\endhead

We  repeat \cite{CM}.
  We   pick $k=1,2,...N$ and
set $f=f_k$ for $k=1$. The other  $f_k$ are defined below. In the
ODE's  there will be error terms    of the form
$$E_{ODE}(k)=   O( |z |^{2 N+2 }
 )+O(  z^{N+1}   f _{k}  )  +O(f^2_{k})+
O(\beta (|f _{k}|^2)f _{k} ).$$ In the PDE's there will be error
terms   of the form
$$E_{PDE}(k)=   O_{loc}( |z |^{N+2}
 )  +O_{loc}(  z f _{k}  )+O_{loc}(f^2_{k})+
O(\beta (|f _{k}|^2)f _{k} ).$$ In the right hand sides of the
equations  (2.3-4) we substitute $\dot \gamma $ and $\dot \omega $
using the modulation equations. We repeat the procedure a sufficient
number of times until we can write for $k=1$ and $ f_1=f$

$$\aligned   i\dot \omega \frac{d\| \phi _\omega \| _2^2}{d\omega  }
=&\langle   \sum _{   m +n =2}^{ 2N+1 } z^m \bar z^n\Lambda
_{m,n}^{(k)}(\omega ) +   \sum _{   m +n =1}^{ N } z^m \bar z^n
A_{m,n}^{ (k) }(\omega ) f_{k}  +E_{ODE}(k) , \Phi (\omega ) \rangle
  \\  i\dot z -\lambda  z =&\langle
\text{ same as above }, \sigma _3 \xi (\omega )  \rangle
\\ i\partial _t f_{k}=&\left ( H _{\omega}  +\sigma _3 \dot \gamma
\right  )f_{k} + E_{PDE}(k)+ \sum _{k+1\le  m +n \le N +1}z^m \bar
z^n R_{m,n}^{(k)}(\omega ),
\endaligned
$$
with  $A_{m,n}^{(k)} $, $R_{m,n}^{(k)}     $  and $\Lambda
_{m,n}^{(k)} (\omega , x) $ real
 exponentially  decreasing  to 0 for $|x|\to \infty$ and continuous in $(\omega , x) $.
Exploiting $|(m-n)\lambda (\omega )|<\omega $ for $  m+n \le N$,
$m\ge 0$, $n\ge 0$,
 we define inductively $f_k$ with $k\le N$ by

 $$f_{k-1}=-\sum _{  m+n=k}z^m \bar z^n R_{H_{\omega }}((m-n)\lambda (\omega ) ) R_{m,n}^{(k-1)}   (\omega
 )
+f_k.$$ Notice that if $ R_{m,n}^{(k-1)}   (\omega
 ,x)$ is real
 exponentially  decreasing  to 0 for $|x|\to \infty$, the same is
 true for $R_{H_{\omega }}((m-n)\lambda (\omega ) ) R_{m,n}^{(k-1)}   (\omega
 )$ by $|(m-n)\lambda (\omega )|<\omega $. By induction $f_k$ solves the above equation with
the above notifications. Now we manipulate the   equation for $f_N$.
We fix $\omega _1=\omega (0)$. We write

$$\aligned &i\partial _t P_c(\omega _1)f_{N}=\left \{ H_{\omega _1}   + (\dot \gamma +\omega -\omega _1)
(P_+(\omega _1)-P_-(\omega _1))\right \} P_c(\omega _1)f_{N} +\\&
+P_c(\omega _1)\widetilde{E}_{PDE}(N) +\sum _{  m+n= N+1} z^m \bar
z^n P_c(\omega _1)R_{m,n}^{(N)} (\omega _1)
\endaligned \tag 4.1$$ where we split $P_c(\omega _1)=P_+(\omega _1)+P_-(\omega
_1)$ with $P_\pm (\omega _1)$ the projections in $\sigma
_c(H_{\omega _1}) \cap \{ \lambda : \pm \lambda \ge \omega _1 \} $,
see \cite{BP2,BS,C2} and Appendix B, and with
$$\aligned & \widetilde{E}_{PDE}(N)= E_{PDE}(N)  +  \sum _{  m+n= N+1} z^m \bar
z^n  \left ( R_{m,n}^{(N)} (\omega  )-R_{m,n}^{(N)} (\omega _1)
\right ) +  \varphi (t,x  ) f_N\\& \varphi (t,x  ):=
 (\dot \gamma +\omega -\omega _1) \left (P_c(\omega _1)\sigma _3-  (P_+(\omega _1)-P_-(\omega _1))  \right ) f_N +
   \left (  V  (\omega ) -
 V  (\omega _1) \right )   f_N
\\& +(\dot \gamma +\omega -\omega _1) \left ( P_c(\omega  )-  P_c(\omega _1)\right )\sigma _3f_N.
  \endaligned \tag 4.2
$$
By Appendix B   for $C_N(\omega _1)$ upper semicontinuous in $\omega
_1$, $\forall $ $N$ we have
$$ \|  \langle x \rangle ^{N}  (P_+(\omega
_1)-P_-(\omega _1)-P_c(\omega _1)\sigma _3) f\|  _{L^2_x}\le
C_N(\omega _1) \|  \langle x \rangle ^{-N}    f\|  _{L^2_x},\tag
4.3$$ see also \cite{BP2,BS}. Then $\varphi (t,x  )$  can be treated
as a small cutoff function.  We write
$$\aligned &f_{N}=-\sum _{  m+n=N+1} z^m \bar z^nR_{H_{\omega _1}}((m-n)  \lambda (\omega _1) +i0  )
 P_c(\omega
_1)R_{m,n}^{(N)}   (\omega _1)  +  f_{N+1}. \endaligned \tag 4.4$$
   Then
$$\aligned &i\partial _t P_c(\omega _1) f_{N+1}=\left ( H_{\omega _1}    + (\dot \gamma +\omega -\omega _1)
(P_+(\omega _1) -P_-(\omega _1) )\right ) P_c(\omega _1) f_{N+1} +
\\& +  \sum _{\pm } O(\epsilon |z|^{N+1}
)R_{H_{\omega _1}} (\pm  (N+1) \lambda (\omega _1)+i0 ) R_\pm
(\omega _1) +P_c(\omega _1)\widehat{ {E}}_{PDE}(N)
\endaligned \tag 4.5$$ with $R_+ =R_{ N+1,0}^{(N)} $ and $R_- =R_{
0, N+1 }^{(N)} $ and $\widehat{ {E}}_{PDE}(N)=\widetilde{
{E}}_{PDE}(N)+ O_{loc}(\epsilon z^{N+1})$, where we have  used that
$( \omega -\omega _1)=O(\epsilon )$ by Theorem 2.1. Notice that
$R_{H_{\omega _1}} (\pm (N+1) \lambda (\omega _1)+i0 ) R_\pm (\omega
_1)\in L^\infty$ do not decay spatially.
  In the ODE's with $k=N$, by the standard theory of normal forms
and following the idea in Proposition 4.1 \cite{BS}, see \cite{CM}
for details, it is possible to introduce new unknowns
$$\aligned &\widetilde{\omega }=\omega +q(\omega , z,\bar z) +\sum _{1\le  m
+n \le N }z^m  \bar z^n\langle f_N,\alpha _{mn}(\omega )\rangle ,\\&
\widetilde{z } =z  +p (\omega , z,\bar z) +\sum _{1\le  m +n \le N
}z^m \bar z^n\langle f_N,\beta _{mn}(\omega )\rangle ,\endaligned
\tag 4.6$$ with $p (\omega , z,\bar z)=\sum p_{ m,n}(\omega )
z^m\bar z^n$ and $q(z,\bar z)=\sum q_{ m,n}(\omega ) z^m\bar z^n$
polynomials in $(z,\bar z)$ with real coefficients and $O(|z|^2)$
near 0, such that we get
$$\aligned & i  \dot {\widetilde{\omega}} =  \langle  {E}_{PDE}(N) , \Phi \rangle   \\ &
 i\dot {\widetilde{z} }-\lambda  (\omega )\widetilde{z}  =
 \sum _{ 1\le  m \le N} a_{ m  }(\omega )|\widetilde{z}^{m}| ^{2 }
 \widetilde{z} +\langle  E_{ODE}(N) , \sigma _3 \xi   \rangle +\\& +
  \overline{\widetilde{z}}^N \langle  A_{0,N}^{(N)}(\omega ) f _{N} , \sigma _3
\xi   \rangle  .
\endaligned \tag 4.7$$
with $   a_{ m  }(\omega )  $ real. Next step is to substitute $
f_{N} $ using (4.4). After eliminating by a new change of variables
$\widetilde{z}= \widehat{{z}}+p(\omega
,\widehat{z},\overline{\widehat{z}})$ the resonant terms, with
$p(\omega ,\widehat{z},\overline{\widehat{z}})=\sum \widehat{p}_{
m,n}(\omega ) z^m\bar z^n$ a polynomial in $(z,\bar z)$ with real
coefficients $O(|z|^2)$ near 0, we get

$$\aligned & i  \dot {\widehat{\omega}} =  \langle  {E}_{PDE}(N) , \Phi \rangle   \\ &
 i\dot {\widehat{z} }-\lambda  (\omega )\widehat{z}  =
 \sum _{ 1\le m\le N} \widehat{a}_{ m  }(\omega )|\widetilde{z}^{m}| ^{2 }
 \widehat{z} +\langle  E_{ODE}(N) , \sigma _3 \xi   \rangle   -\\&
 -
  |\widehat{z} ^N|^2\widehat{z}  \langle  \widehat{A}_{0,N}^{(N)}(\omega )
  R_{H_{\omega _1}}((N+1) \lambda  (\omega _1) +i0)P_c(\omega
_0)R_{N+1 ,0}^{(N)}(\omega _1)     ,\sigma _3\xi   \rangle
\\& +  \overline{\widehat{z}} ^N \langle  \widehat{A}_{0,N}^{(N)}(\omega ) f
_{N+1} , \sigma _3 \xi   \rangle
\endaligned \tag 4.8$$
with $\widehat{a}_{ m  }$, $\widehat{A}_{0,N}^{(N)}$ and
$R_{N+1,0}^{(N)}$ real. By $\frac{1}{x-i0}=PV\frac{1}{x}+i\pi \delta
_0(x)$ and by \cite{BP2,BS} we can denote by $\Gamma  (\omega ,
\omega _0)$ the quantity
$$\aligned &\Gamma  (\omega ,
\omega _1)= \Im \left (\langle  \widehat{A}_{0,N}^{(N)}(\omega )
  R_{H_{\omega  _1}}((N+1) \lambda (\omega _1) +i0)P_c(\omega
_1)R_{N+1,0}^{(N)}(\omega _1 \sigma _3 \xi   (\omega )\rangle \right
)\\& =\pi
  \langle  \widehat{A}_{0,N}^{(N)}(\omega )
  \delta (H_{\omega  _1}- (N+1) \lambda (\omega _1)  )P_c(\omega
_1)R_{N+1,0}^{(N)}(\omega _1) \sigma _3 \xi   (\omega )\rangle .
 \endaligned $$
Now we assume the following: \proclaim{Hypothesis 4.2} There is a
fixed constant $\Gamma >0$ such that $|\Gamma  (\omega , \omega
)|>\Gamma .$
\endproclaim
Notice that the FGR hypothesis in \cite{GS1} asks $\Gamma (\omega ,
\omega  ) >0.$ We will prove in Corollary 4.7  that in fact $ \Gamma
(\omega , \omega ) >\Gamma  $. By continuity and by Hypothesis 4.2
we can assume $|\Gamma (\omega , \omega _1 )|>\Gamma /2.$ Then we
write

$$\aligned &  \frac{d}{dt} \frac{|{\widehat{z} }|^2}{2} =    -\Gamma  (\omega , \omega  _1)
  |z|  ^{2N+2} +\Im \left (
\langle  \widehat{A}_{0,N}^{(N)}(\omega )f _{N+1} , \sigma _3 \xi
(\omega )\rangle \overline{\widehat{z}} ^{N+1}   \right ) \\&+ \Im
\left ( \langle E_{ODE}(N) , \sigma _3 \xi  (\omega ) \rangle
\overline{\widehat{z}} \right ) .
\endaligned \tag 4.9
$$

\head \S 4.2 Nonlinear estimates \endhead

By an elementary continuation argument, the following a priori
estimates imply  inequality (1) in Theorem 4.1, so to prove (1) we
focus on:

\proclaim{Lemma 4.3} There are fixed  constants $C_0$ and $C_1$ and
$\epsilon _0>0$  such that   for any $0<\epsilon \le \epsilon _0$
  if we have
$$\| \widehat{ z} \| _{L^{2N+2 }_t}^{N+1} \le 2C_0\epsilon \quad \&
 \quad \| f_N \| _{L^\infty  _tH ^{1 }_x \cap L^5 _tW ^{1,10}_x
 \cap L^4_tL^\infty _x\cap L^2_tH^{1,-2}_x } \le 2C_1\epsilon  \tag 4.10$$ then we
obtain the improved inequalities
$$  \align&
\| f_N \| _{L^\infty  _tH ^{1 }_x \cap L^5 _tW ^{1,10}_x \cap
L^4_tL^\infty _x \cap  L^2_tH^{1,-2}_x} \le C_1\epsilon  , \tag 4.11
\\& \| \widehat{z} \| _{L^{2N+2 }_t}^{N+1} \le C_0\epsilon .\tag
4.12 \endalign$$
\endproclaim
{\it Proof}. Set $\ell (t):=\gamma +\omega -\omega _0$.  First of
all, we have:

\proclaim{Lemma 4.4} Let $g (0,x)\in H^1_x\cap L^2_c(\omega _1)$ and
let  $\omega (t)$ be  a continuous function. Consider    $  i g _t=
\left \{ H _{\omega _1}   + \ell (t) (P_+(\omega _1)-P_-(\omega
_1))\right \}g+P_c(\omega _1) F.$  Then for a fixed $C=C(\omega _1)$
upper semicontinuous in $\omega _1$ we have    $$ \| g \| _{L^\infty
_tH ^{1 }_x \cap L^5 _tW ^{1,10}_x \cap L^4_tL^\infty _x}
 \le  C   \| g(0,x)\| _{H^1}+C \| F\|  _{   L^1
_tH^1 _x  + L ^{\frac{4}{3}} _tW ^{1,1} _x + L^2_tH^{1, 2}_x} .$$
\endproclaim
Lemma 4.4 follows from Lemmas 3.1 and 3.5 and $ P_\pm (\omega
_1)g(t)=$
$$=e^{-it H_{\omega _1}}e^{- i\int _0^t \ell  (\tau ) d\tau }P_\pm
(\omega _1)g(0)  -i\int _0^{t }e^{-i(t-s) H_{\omega _1}} e^{\pm
i\int _s^t \ell (\tau ) d\tau }P_{ \pm}(\omega _1)
  F  (s) ds
$$

 \proclaim{Lemma 4.5} Consider equation (4.1) for
 $f_N$ and assume (4.10). Then we can split $\widetilde{E}_{PDE}(N)=
 X+ O(f_N^5)$ such that
 $ \|  X\| _{   H_x^{1,M} L_t^2 }  \lesssim \epsilon ^2
 $ for any fixed $M$ and  $ \|O(f_N^5)\| _{ L^1
_tH^1 _x}   \lesssim \epsilon ^5.
 $
\endproclaim
{\it Proof of Lemma 4.5.} Schematically we have for a cutoff $\psi
(x)$
$$\widetilde{E}_{PDE}(N) = O(\epsilon ) \psi (x) f_N +O_{loc}( |z |^{N+2}
 )  +    O_{loc}(  z f _{N}  )+O_{loc}(f^2_{N})+
O(\beta (|f _{N}|^2)f _{N} ).$$ By (4.10)  for all the terms in
$\widetilde{E}_{PDE}(N)$ except the last one  and whose sum we call
$X$, we have:\medskip {\item{(1)}} $ \| \la x \ra ^{ M}O(\epsilon )
\psi (x) f_N \| _{ H_x^1 L_t^2 } \lesssim \epsilon \| \la x \ra ^{
-5} f_N \| _{H_x^1 L_t^2}\lesssim \epsilon ^2;$
\medskip {\item{(2)}}$ \| \la x \ra ^{
M} O_{loc}(  z f _{N}  ) \| _{ H_x^1 L_t^2 } \lesssim \| z\|
_{\infty} \| \la x \ra ^{ -5} f_N \| _{H_x^1 L_t^2}\lesssim \epsilon
^2 ;$
\medskip
{\item{(3)}} $ \| \la x \ra ^{ M}O_{loc}(    f _{N}^2  ) \| _{ H_x^1
L_t^2 } \lesssim   \| \la x \ra ^{ -5} f_N \| _{H_x^1
L_t^2}^2\lesssim \epsilon ^2.$ \medskip {\item{(4)}} $ \| \la x \ra
^{ M}O_{loc}(    |z |^{N+2} ) \| _{ H_x^1 L_t^2 } \lesssim
   \epsilon  \|z  ^{N+1}\| _{L^2_t} \lesssim \epsilon ^2 .$

  \noindent This yields $ \| \la x \ra ^{  M} X\| _{   H_x^1 L_t^2 }  \lesssim \epsilon ^2
 $.  Observe that schematically $\| \beta (|f _{N}|^2)f _{N}\|
_{W^{1,r}_x}\lesssim  \| f_N^5\| _{W^{1,r}_x}$ for all $ r\in
(1,\infty)$, if on the right hand side we mean all the fifth powers
of the components of $f_N$. Then we have

  $$ \|      f _{N}^5 \| _{  L_t^1H^1_x } \lesssim \left \|
\| f_N\| _{W^{1,10}_x}    \| f_N\| ^4_{L^{ 10}_x} \right \| _{L^1_t
} \le \| f_N\| _{L^5_tW^{1,10}_x}^5 \lesssim \epsilon ^5.\tag 5$$

\bigskip

 {\it Proof of (4.11).} Recall that $f_N$ satisfies equation (4.1)
whose right hand side is $P_c(\omega _1)\widetilde{E}_{PDE}(N) +
O_{loc}(    z  ^{N+1}  )$. In addition to Lemma 4.5 we have the
estimate $ \| O_{loc}(    z  ^{N+1}  ) \| _{ L_t^2H_x^{1,M} }
\lesssim \| z \| _{ L_t^{2N+1} }^{N+1} \lesssim 2C_0\epsilon  .$ So
by Lemmas 3.1-4, for some fixed $c_2$ we get schematically
$$ \| f_N \| _{L^\infty  _tH ^{1 }_x \cap L^5 _tW ^{1,10}_x \cap L^4_tL^\infty
_x } \le 2c_2C_0\epsilon +\| f_N(0) \| _{ H ^{1 }_x}  +O(\epsilon
^2)$$ where $ \| f_N(0) \| _{ H ^{1 }_x} \le c_2 \epsilon$ for fixed
$c_2\ge 1$, $O(\epsilon ^2)$ comes from all the  terms on the right
of (4.1) save for
 the $ R_{m,n}^{(N)}   (\omega _0) z^m  \bar z^n $ terms which contribute the  $2c_2C_0\epsilon
 $.
Let now $f_N= g+h$ with
$$\aligned & i g _t=
\left \{ H _{\omega _1}   + \ell (t) (P_+(\omega _1)-P_-(\omega
_1))\right \}g+X \, , \quad  g(0)=f_N(0)\\&   i h _t= \left \{ H
_{\omega _1}   + \ell (t) (P_+(\omega _1)-P_-(\omega _1))\right
\}h+O(f_N^5) \, , \quad  h(0)=0\endaligned $$ in the notation of
Lemma 4.5. Then   $\|  g \| _{H ^{1,-2} _xL^2_t} \lesssim
2C_0\epsilon +O(\epsilon ^2)+c_0\epsilon $  by Lemmas 3.3-4 for a
fixed $c_0$. Finally by Lemma 3.3
$$\aligned &\int _0^\infty  \|     e^{-i(t-s) H_{\omega _1}}  e^{\pm i \int _s^t\ell (\tau ) d\tau }
O(f_N^5)  (s)\| _{H^{1,-2}_xL^2_t}
 \lesssim \int _0^\infty \|
O(f_N^5) (s) \| _{ H^{1}_x}\lesssim \epsilon ^5.
\endaligned $$
 So if we set  $C_1\approx 2C_0+1$ we obtain (4.11). We
need to bound $C_0$.

\medskip

{\it Proof of (4.12).}
 We first need: \proclaim{Lemma 4.6} We can
decompose $f _{N+1}= h_1+h_2+h_3 +h_4$ with   for  a fixed large
$M>0$: {\item {(1)}}       $\| \langle x\rangle ^{-M} h _1\|
_{L^2_{tx}} \le O(\epsilon ^{2})   ;$ {\item {(2)}}   $\| \langle
x\rangle ^{-M} h _2\| _{L^2_{tx}}\le O(\epsilon ^{2}) ;$

{\item {(3)}}   $\| \langle x\rangle ^{-M} h _3\| _{L^2_{tx}} \le
O(\epsilon ^{2});$  {\item {(4)}}   $\| \langle x\rangle ^{-M} h
_4\| _{L^2_{tx}} \le  c(\omega _1)\epsilon  $ for a fixed $c(\omega
_1) $ upper semicontinuous in $\omega _1$.
\endproclaim
{\it Proof of Lemma 4.6}. We set $$\aligned & i\partial _t h_1=\left
( H_{\omega _1}   + \ell (t) (P_+-P_-)\right ) h_1\\& h_1(0)=\sum _{
m+n=N+1} R_{H_{\omega _1}}((m-n)\lambda (\omega _1)+i0)
R_{m,n}^{(N)} (\omega _1)z^m(0)  \bar z^n (0) .\endaligned $$  We
get $\| \langle x\rangle ^{-M} h _1\| _{L^2_{tx}} \le c(\omega _1)
|z(0)|^2 \sum \| \langle x \rangle ^{\gamma } R_{m,n}^{(N)} (\omega
_1)\| _{L^2_x}=O(\epsilon ^2) $ by the inequality (4.13) below, see
\cite{BP1,BS}, which says that for any $\gamma > \gamma _0$ for some
given $\gamma _0$,

$$\aligned &
\| \langle x \rangle ^{-\gamma } e^{-iH _\omega t}R_{H_\omega}
(\Lambda +i0) P_c(\omega ) g\| _2 < C(\Lambda , \omega )\langle t
\rangle ^{-\frac 32} \| \langle x \rangle ^{\gamma } g \| _2 \, ,\,
\Lambda >\omega   , \endaligned \tag 4.13
$$
with   $ C(\Lambda , \omega )$ upper semicontinuous in $\omega $ and
in $\Lambda  $. Next, we set $h_2(0)=0$ and

$$\aligned &i\partial _t h_{2}=\left ( H_{\omega _1}   + \ell (t)  (P_+-P_-)\right ) h_{2}
+\\& + O(\epsilon z^{N+1}) R_{H_{\omega _1}}((N+1)\lambda (\omega
_1)+i0 ) R_{N+1,0}^{(N)} (\omega   _0)
\\& +
O(\epsilon z^{N+1}) R_{H_{\omega _1}}(-(N+1)\lambda (\omega _1)+i0 )
R_{0, N+1}^{(N)} (\omega   _1) .\endaligned $$ Then we have
$h_{2}=h_{21 }+h_{22}$ with $h_{2j}=\sum _\pm h_{2j\pm }$ with $h_{2
1\pm } (t)=$

$$\int _0^te^{-iH_{\omega _1}(t-s)}e^{\pm i\int _s^t \ell(\tau )
d\tau }P_{ \pm}
 z^{N+2}(s) R_{H_{\omega _1}}((N+1)\lambda (\omega _1)+i0 )  R_{N+1,0}^{(N)}   (\omega   _1) ds $$
and  $h_{22\pm }$ defined similarly but with $R_{H_{\omega
_1}}(-(N+1)\lambda (\omega _1)+i0 ) R_{0,N+1 }^{(N)}$ . Now by
(4.13) we get
$$\|   \langle x \rangle ^{-M } h_{2j\pm } (t) \| _{L^2_{x} }  \le C \epsilon \int _0^t \langle
t-s\rangle ^{-\frac 32} |z(s)|^{N+1}  ds$$ and so $  \| \langle
x\rangle ^{-M} h _2\| _{L^2_{tx}} \le \epsilon \| z \|  _{L^{2N+2}
_{t }} ^{N+2}= O(\epsilon ^3).$ Let   $h_3(0) =0$   and

$$\aligned &i\partial _t P_c(\omega _1) h_3=\left ( H_{\omega _1}   + \ell (t)
(P_+(\omega _1) -P_-(\omega _1) )\right ) P_c(\omega _1) h_3 +
 P_c(\omega _1) \widetilde{E}_{PDE}(N). \endaligned
 $$
Then by the argument in the proof of (4.11) we get claim (3).
Finally let   $h_4(0) =f_N(0)$   and

$$\aligned &i\partial _t P_c(\omega _1) h_4=\left ( H_{\omega _1}   + \ell (t)
(P_+(\omega _1) -P_-(\omega _1) )\right ) P_c(\omega _1) h_4  .
\endaligned
 $$
Then by  Lemma 3.3   $ \| \langle x\rangle ^{-M} h _4\| _{L^2_{tx}}
 \lesssim \| f_N(0)\| _{L^2_x}\le c(\omega _1)\epsilon $ we get
(4).
\bigskip
{\it Continuation of proof of Lemma 4.3}. We integrate (4.9) in
time. Then by Theorem 2.1  and by Lemma 4.4 we get, for $A_0$   an
upper bound  of the    constants $A_0(\omega )$ of Theorem 2.1,

$$ \| \widehat{z}\| _{L^{2N+2}_t}^{2N+2}\le  A_0\epsilon ^2+2c(\omega _1)\epsilon
\| \widehat{z}\| _{L^{2N+2}_t}^{N+1} + o(\epsilon ^2).$$ Then we can
pick    $C_0=( A_0+2c(\omega _1 +1)$ and this proves that (4.10)
implies (4.12). Furthermore $\widehat{z}(t)\to 0$ by
$\frac{d}{dt}\widehat{z}(t)=O(\epsilon ).$

\bigskip

As in \cite{CM} in the above argument we did not use the sign of
$\Gamma (\omega , \omega _1)$. As in \cite{CM}  it is nonnegative.

 \proclaim{Corollary 4.7} If Hypothesis 4.2 holds, then $\Gamma (\omega , \omega
 )>\Gamma $.
 \endproclaim Suppose we have $\Gamma (\omega , \omega _1
 )<-\Gamma $. We can pick initial datum so that $f_{N+1}(0)=0$ and $z(0)\approx \epsilon $. Then
 following the proof of Lemma 4.6, by integrating (4.9)   and using $h_4=0$,    we   get
 $$|\widehat{z}(t)|^2-|\widehat{z}(0)|^2\ge \Gamma \int _0^t|\widehat{z}|^{2N+2}+ o(\epsilon ) \left (  \int
 _0^t|\widehat{z}|^{2N+2} \right ) ^{\frac{1}{2}}+ o(\epsilon ^2).$$
For large  $t$ we have
 $|\widehat{z}(t)|<|\widehat{z}(0)|$ since $ z(t) \to 0$, so
for large $t$ we get $\int _0^t|\widehat{z}|^{2N+2}=o(\epsilon
 ^2).$ In particular for $t\to \infty $ we get
$   \epsilon ^2\le   o(\epsilon ^2) $ which is absurd for $\epsilon
\to 0 $.

\bigskip
The proof that, for $^tf_N(t)=(h(t),\overline{h}(t))$, $h(t)$ is
asymptotically free for $t\to \infty$, is similar to the analogous
one in \cite{CM} and we skip it.

\head \S Appendix A.  Orbital stability: sketch of proof of Theorem
2.1\endhead

We sketch the proof of Theorem 2.1.
 \proclaim{Lemma A.1}  Suppose that $e^{i\omega t} \phi
_ {\omega } (x)$ satisfies (H6). In dimension $n>1$ assume also that
$$L_+=-\Delta +\omega - \beta (  \phi _\omega ^2  )-2\beta ' ( \phi
_\omega ^2 )\phi _\omega \tag A.1$$ has exactly one negative
eigenvalue. Then $\exists \, \epsilon >0$ and a $A_0(\omega )>0$
such that for any  $\| u(0,x) - \phi _ {\omega } \| _{H^1(\Bbb
R^n)}<\epsilon $ we have for the corresponding solution
$$\inf  \{  \|  u(t,x) -e^{i \gamma }\phi _ {\omega } (x-x_0) \| _{H^1(x\in \Bbb R^n)}
: \gamma \in \Bbb R \, \& \, x_0 \, \in \Bbb R^n \}  < A_0(\omega )
\epsilon .
$$\endproclaim

The proof consists in   the argument in \cite{We1} with a minor
change due to D. Stuart \cite{S}. We have invariants:
$$ \aligned & Q(f)= \frac{1}{2} \int  _{\Bbb R^{n}} |f(x)|^2 dx \, , \, M(f)= \frac{1}{2} \Im
\int _{\Bbb R^{n}}  \overline{f(x)} \nabla f  (x)  dx \, , \\& E(f)
=\int _{\Bbb R^{n}} \left ( \frac{| \nabla f (x) |^2}{2}+ F(|f|)
\right ) dx.\endaligned$$ For $\Theta (t,x)= \dfrac{v\cdot x }{2} +
\vartheta (t)$ we have

$$\aligned &  M(e^{i\Theta }f)= \frac{1}{2} \Im \int  _{\Bbb R^{n}}  e^{- i\Theta } \overline{f(x)} e^{i\Theta }(  \nabla f (x)  +i
\frac{v}{2}f(x)) dx= M(f)+ \frac v2 Q(f) \\& E(e^{i\Theta }f) =\int
_{\Bbb R^{n}} \left ( \frac{| \nabla f (x)+ i \frac{v}{2}f(x)
|^2}{2}+ F(|f|) \right ) dx=E(f)+ \frac{v^2}{4}Q(f)+v\cdot M(f).
\endaligned $$
We define now from the invariants of motion
$$H(u)=E(u)+\omega (t) Q(u)-v(0)\cdot M(u)= E(u_0)+\omega (t)
Q(u_0)-v(0)\cdot M(u_0)$$ with $   v(0)$ initial velocity, $
 \omega (t) $ a function defined later, $
 u_0(x)=u(0,x)  .  $ The idea of choosing $   v(0)$ is in
 \cite{S}.
For $ y$ the coordinate in the moving frame, we consider the ansatz
$ u=e^{i\Theta }( \phi _\mu (y) +r(t,y))$  satisfying the usual
modulation equations
$$\langle Q'  ( \varphi _\mu ) , r(t)\rangle = \langle M' ( \varphi _\mu ) , r(t)\rangle
=0.$$    After the above preparation  we start the usual expansion
$$\aligned & H(e^{i\Theta }( \phi _\mu   +r)) = E(e^{i\Theta }( \phi _\mu
+r)) + \omega Q(\phi _\mu   +r)-v(0)\cdot  M(e^{i\Theta }( \phi _\mu
+r)) =\\& = E(\phi _\mu +r) + \left ( \omega  + \frac{v^2-2v(0)\cdot
v}{4} \right ) Q(\phi _\mu +r)+(v-v(0)) \cdot M( \phi _\mu +r)=\\& =
E(\phi _\mu +r) + \left ( \omega - \frac{v^2(0)}{4}  +
\frac{(v-v(0))^2}{4} \right ) Q(\phi _\mu +r)+(v-v(0)) \cdot M( \phi
_\mu +r).
\endaligned $$
Define $\omega =\dfrac{v^2(0)}{4} +\mu .$  Then, setting $d(\mu) =E
(\phi _\mu
 )+\mu  Q(\phi _\mu
 ) $ and   $q(\mu )
=Q(\phi _\mu
 )$
$$\aligned & H(u)=d(\mu ) + \frac{(v-v(0))^2}{4}   q(\mu )+\langle E' (\phi
_\mu ) + \mu Q' (\phi _\mu )+(v-v(0)) M' ( \phi _\mu ),r \rangle
+\\& + \frac 12 \langle \left [ E'' (\phi _\mu ) + \mu Q'' (\phi
_\mu ) \right ] r,r \rangle + \frac{(v-v(0))^2}{4} \langle
  Q' (\phi _\mu )+ \frac{ Q''}{2} (\phi _\mu )  r ,r \rangle +\\& +\frac{1}{2}          (v-v(0))
 \cdot  \langle M'' (
\phi _\mu )r,r \rangle +o(\| r\| _{H^1}^2).
\endaligned$$
From modulation and from $E' (\phi _\mu ) + \mu Q'(\phi _\mu )=0$ we
get
$$\langle E' (\phi
_\mu ) + \mu Q'(\phi _\mu )+(v-v(0)) \cdot M'( \phi _\mu ),r \rangle
=0.$$ So
$$H(u)=d(\mu ) + \frac{(v-v(0))^2}{4}   q(\mu )+ \frac 12 \langle \left [ E'' (\phi _\mu ) + \mu
Q'' (\phi _\mu ) \right ] r,r \rangle + o(\| r\| _{H^1}^2).$$
Proceeding similarly
$$\aligned &E(u_0)+\omega (t)
Q(u_0)-v(0)\cdot M(u_0)  =d(\mu (0)) + \\&   \frac 12 \langle \left
[ E'' (\phi _{\mu (0)}) +            {\mu (0)}  Q'' (\phi _{\mu (0)}
) \right ] r(0),r (0)\rangle + (\mu -{\mu (0)}) q ( {\mu (0)} )
+o(\| r(0)\| _{H^1}^2) .\endaligned$$  Recall now that $d' ( {\mu
(0)} )=q ( {\mu (0)} )$ so by equating the last two displayed
formulas and after Taylor expansion of $d(\mu )$ we get the
following result:
$$\aligned &\frac{d'' (\mu  (0))}{2} (\mu -\mu (0))^2 + \frac{(v-v(0))^2}{4}   q(\mu )
+ \frac 12 \langle \left [ E'' (\phi _\mu ) + \mu Q'' (\phi _\mu )
\right ] r,r \rangle \le \\& \le \frac 12 \langle \left [ E'' (\phi
_{\mu (0)}) +            {\mu (0)} Q'' (\phi _{\mu (0)} ) \right ]
r(0),r(0)\rangle +o(\| r\| _{H^1}^2) +o(\| r(0)\| _{H^1}^2) .
\endaligned
$$
This implies $(\mu -\mu (0))^2 +  (v-v(0))^2 +\| r\| _{H^1}^2 \le C
\| r(0)\| _{H^1}^2 $ because of the fact that
$$ \langle \left [ E'' (\phi _\mu ) + \mu Q'' (\phi _\mu )
\right ] r,r \rangle \approx \| r  \| _{H^1}^2.$$

\head \S Appendix B.   Proof of estimate (4.3)\endhead

\proclaim{Lemma B.1} The following operators  $P_\pm (\omega )$ are
well defined:
$$\aligned &P_+(\omega )u =\lim _{\epsilon \to 0^+}
 \frac 1{2\pi i}
\lim _{M \to +\infty} \int _\omega ^M \left [ R(\lambda +i\epsilon
)- R(\lambda -i\epsilon )
 \right ] ud\lambda \\&
P_-(\omega )u =\lim _{\epsilon \to 0^+}
 \frac 1{2\pi i}
\lim _{M \to +\infty} \int _{-M }^{-\omega } \left [ R(\lambda
+i\epsilon )- R(\lambda -i\epsilon ) \right ] ud\lambda .
\endaligned
$$
  For any
$M>0$ and $N>0$ and for $C=C (N,M,\omega  )$ upper semicontinuous in
$\omega $, we have $$  \|  \langle x \rangle ^{M}  (P_+(\omega
 )-P_-(\omega  )-P_c(\omega  )\sigma _3) f\|  _{L^2 }\le
C  \|  \langle x \rangle ^{-N}    f\|  _{L^2 }.\tag 2$$
\endproclaim
(2) for $M= 2$      is stated in \cite{BP2} with the proof sketched
in \cite{BS}.

{\it Proof.} The first part is a consequence of \cite{KS}. We prove
(2) following the argument in \S 7 \cite{C2}. For this proof we set
$ L^2_s =\langle x\rangle ^{-s}L^2$, $H=H_\omega$, $H_0=\sigma _3
(-\Delta + \omega )$,
  $R_0(z)= (H_0-z)^{-1}$ and $R(z)= (H-z)^{-1}$.
To prove (2) it is enough to write $P_c=P_++P_-$ and to prove
$\|\left [ P_\pm \sigma _3 \mp P_\pm \right ] g\|  _{L^2_{ M}}  \le
c\| g\| _{L^2_{ -N}} .$ It is not restrictive to consider only
$P_+$. Setting $H=H_0+V$, we write
$$\aligned &\sum_\pm \pm R(\lambda \pm  i\epsilon )
= \sum _\pm \pm (1+ R_0(\lambda \pm  i\epsilon ) V )^{-1}
R_0(\lambda \pm  i\epsilon ) .
\endaligned \tag 3
$$
By elementary computation
$$  R_0(\lambda \pm  i\epsilon )
\sigma _3 =R_0(\lambda \pm  i\epsilon ) -2 (-\Delta +\omega +\lambda
\pm i \epsilon )^{-1} \text{diag}(0,1).$$ Therefore  $$ \text{rhs}
\,(3) \sigma _3 = \text{rhs} \,(3)+2\sum _\pm \pm (1+ R_0(\lambda
\pm i\epsilon ) V )^{-1} \text{diag}(0,1)
 (-\Delta +\omega +\lambda \pm i
\epsilon )^{-1}.
$$
Hence we are reduced to show that
$$Ku=
\lim _{\epsilon \to 0^+} \lim _{M \to +\infty} \sum _\pm \pm \int
_\omega ^M (1+ R_0(\lambda \pm  i\epsilon ) V )^{-1}
\text{diag}(0,1)
 (-\Delta +\omega +\lambda \pm i
\epsilon )^{-1}
 u d\lambda
$$
defines an operator such that  for some fixed $c$  $$\| Ku \|
_{L^2_{ M}} \le c \| u \| _{L^2_{ -N}} \tag 4$$
 For $m \ge 1$ we expand $(1+R_0V )^{-1}=\sum
_{j=0}^{m+1} \left [ - R_0V \right ] ^j +R_0VRV (-R_0V)^N $ and we
consider the corresponding decomposition $$K=\sum
_{j=0}^{m+1}K_j^0+\Cal K.\tag 5$$   We have $K_0^0= 0$ since for any
$u\in L^2$ we have
$$  \lim _{\epsilon \to 0^+} \lim _{M \to +\infty} \int
_\omega ^M
 \sum _\pm \pm (-\Delta +\omega +\lambda \pm i
\epsilon )^{-1} \text{diag}(0,1)ud\lambda =0.$$ We next consider
$K_1^0$ and prove $$\| K_1^0u \| _{L^2_{ M}} \le c \| u \| _{L^2_{
-N}}. \tag 6$$ The
 operator $(-\Delta +\omega +z )^{-1}$ has symbol satisfying,
for $\Re z \ge 0 $:
$$\big |\partial _{z }^\beta \partial _{\xi }^\alpha
(|\xi |^2+\omega + z)^{-1} \big | \le c_{\alpha  ,\beta} (|\xi
|+1)^{-|\alpha |} \langle z \rangle ^{-1-\beta }. \tag 7$$ Therefore
we have

$$ \| \langle x \rangle ^{M}(-\Delta +\omega +z )^{-1}u\|_{ L^2}=
 \| \langle  \sqrt{-\Delta _{\xi}} \rangle ^{M}(\xi ^2 +\omega +z )^{-1}\widehat{u}\|_{ L^2}
\le C \langle z \rangle ^{-1 } \| u \|_{ L^2_M}$$ and so for any
$M\in \Bbb R$
$$ \| (-\Delta +\omega +z )^{-1}:L^2_M \to L^2_M\|
\le C \langle z \rangle ^{-1 }.\tag 8$$ We can assume $u$ smooth and
rapidly decreasing. Since for $s>1$ we have $\| R_0(\lambda \pm  i
\epsilon ) :L^2_{s}\to L^2_{-s}\| \le C \langle \lambda \rangle
^{-\frac{1}{2}},$   the following limit is well defined
$$\aligned &
K_1^0u= \lim _{\epsilon \to 0^+} \lim _{M \to +\infty}
 \int _\omega ^M \sum _\pm \pm  \left [ R_0 (\lambda \pm  i
\epsilon )  V   (-\Delta +\omega +\lambda \pm i\epsilon  )^{-1}
\right ] \text{diag}(0,1)
 u d\lambda \\& =
 \int _\omega ^{+\infty}  \left [  R_0 (\lambda + i
0 ) - R_0 (\lambda - i 0 ) \right ] V   (-\Delta +\omega +\lambda
)^{-1} \text{diag}(0,1)
 u d\lambda .\endaligned
$$ By
$
 R_0 (\lambda + i
0 ) - R_0 (\lambda - i 0 )  = 2i\pi \delta (\Delta -\omega +\lambda
) \text{diag}(1,0) $   and for $^tu=(u_1,u_2)$
$$K_1^0u= \int _\omega ^{+\infty }\delta (\Delta -\omega +\lambda  )\text{diag}(1,0) V
(-\Delta +\omega +\lambda  )^{-1} u_2\overrightarrow{e}_2d\lambda .
$$ Up to a constant factor, this is schematically
$$\int _{\Bbb R^2} e^{ix\cdot \xi }
\frac {\hat V( \xi -\eta )\hat u(   \eta )}{\xi ^2+  \eta  ^2
+2\omega} d\eta d\xi . $$ By the correspondence $\partial
_x\leftrightarrow i\xi$ and by Parseval equality, (6) will follow by

$$\aligned & \left \| \int _{\Bbb R }d\eta   \hat u(   \eta ) \hat V^{(\ell _1)}( \xi -\eta  ) \partial _\xi ^{
\ell _2} ({\xi ^2+  \eta  ^2 +2\omega} )^{-1}\right  \| _{L^2_\xi}
\le C(\ell _1,\ell _2)  \|  \hat u  \| _2 \endaligned
$$ which is a consequence of Young inequality. We consider now

$$K_j^0u=(-)^j\lim _{\epsilon \to 0^+} \sum _\pm \pm \int _\omega ^{+\infty }
\left [ R_0(\lambda \pm i \epsilon ) V \right ] ^j
\text{diag}(0,1)(-\Delta +\omega + \lambda \pm i \epsilon ) ^{-1}
ud\lambda .$$   For some $\delta >0$ small but fixed we can deform
the path of integration and write

$$K_j^0u=(-)^j\int _{\omega -\delta -i\infty }^{\omega -\delta +i\infty }
\left [ R_0(\zeta ) V \right ] ^j \text{diag}(0,1)(-\Delta +\omega +
\zeta ) ^{-1}ud\zeta .$$ By (8) we conclude $$\|  K_j^0 u \|
_{L^2_M} \le c \| u \| _{L^2_{-N}}.\tag 9$$ Next we consider also
the reminder term in (5). Arguing as above

$$\aligned &(-)^{m+2}\Cal  K u=
\lim _{\epsilon \to 0^+} \sum _\pm \\&  \pm \int _\omega ^{+\infty }
 R_0(\lambda \pm i \epsilon ) V R(\lambda \pm i \epsilon )
V \left [ R_0(\lambda \pm i \epsilon ) V\right ] ^m
\text{diag}(0,1)(-\Delta +\omega + \lambda \pm i \epsilon ) ^{-1}
ud\lambda \\ &=\int _{\omega -\delta -i\infty }^{\omega -\delta
+i\infty }
 R_0(\zeta ) V  R(\zeta )
V\left [ R_0(\zeta ) V\right ] ^m \text{diag}(0,1)(-\Delta +\omega +
\zeta ) ^{-1}ud\zeta .
\endaligned $$
 For $\Re \zeta =\omega -\delta $, (7) implies $
  (1+|\zeta |)^{-1} \gtrsim \|   R_0(\zeta ) V \colon
L^2_{M} \to L^2_{-N} \| +$
$$ \aligned &
 + \|  V\left [ R_0(\zeta )
V\right ] ^m \colon L^2_{-N} \to L^2_{-N} \| + \| (-\Delta +\omega +
\zeta ) ^{-1} \colon L^2_{-N} \to L^2_{-N} \| .
\endaligned $$ So $ \|  \Cal K
\colon L^2_{M} \to L^2_{-N} \| < \infty $  and this with (6) and (9)
yields (4) and proves $\|\left [ P_+ \sigma _3 - P_+ \right ] u\|
_{L^2_{ M}}  \le c\| u\| _{L^2_{ -N}} .$

\enddocument